\documentclass[Journal]{IEEEtran}

\usepackage{graphicx}
\usepackage{url}
 \usepackage[pdftex]{color}
 \usepackage{mathrsfs}
\usepackage{amsmath}
\usepackage{amsfonts} 
\usepackage{array}
\usepackage{cancel}
\newtheorem{teor}{Theorem}[section]
\newtheorem{corr}{Corollary}[section]

\newtheorem{lemm}{Lemma}[section]

\newtheorem{rem}{Remark}
\usepackage{soul}
\usepackage{tabularx}

\def\bsmat{\left[ \begin{smallmatrix}}
\def\esmat{\end{smallmatrix} \right]}

\newcommand{\tp}{^{\top}}
\newcommand{\pinv}{^{+}}

\newcommand{\beq}{\begin{equation}}
\newcommand{\eeq}{\end{equation}}
\newcommand{\bea}{\begin{eqnarray}}
\newcommand{\eea}{\end{eqnarray}}

\newcommand{\bsea}{\begin{subeqnarray}}
\newcommand{\esea}{\end{subeqnarray}}
\newcommand{\nn}{\nonumber}

\newcommand{\proof}{\noindent {\it Proof. }}

\newcommand{\qed}{\hfill $\Box$ \vskip 2ex}

\def\bmat{\left[ \begin{array}}
\def\emat{\end{array} \right]}

\definecolor{Royalblue}{cmyk}{1,0.30,0.2,0.2}

\definecolor{USred}{rgb}{0.74,0.1,0.1}

\begin{document}

\title{Families of solutions of  algebraic Riccati equations}

\author{Daniele~Alpago, Augusto~Ferrante \thanks{D. Alpago, and A. Ferrante are with   the Department of Information
Engineering, University of Padova, Padova, Italy; e-mail: {\tt\small daniele.alpago@phd.unipd.it} (D. Alpago); {\tt\small augusto@dei.unipd.it} (A. Ferrante).}}

\markboth{DRAFT}{Shell \MakeLowercase{\textit{et al.}}: Bare Demo of IEEEtran.cls for Journals}

\maketitle

\begin{abstract}
We consider Homogeneous Algebraic Riccati Equations  in the general situation when the matrix of the dynamics can be ``mixed''.
We show that in this case the equation may have infinitely many families of solutions.
An analysis of these families is carried over and explicit formulas are derived.
We also derive sufficient conditions under which the union of the families  covers the  whole set of solutions.
\end{abstract}


\section{Introduction}
Since the seminal paper \cite{Willems-71} of J. C. Willems where it was first hinted that (under reasonable assumptions) once given a reference solution of an Algebraic Riccati Equation (ARE), all the others can be parametrized in terms of the solutions of an associated Homogeneous Algebraic Riccati Equation (HARE), a huge amount of literature has been produced on this topic.
In particular, \cite{Coppel,Shayman} developed the work of Willems establishing what is referred to as the Willems-Coppel-Shayman parametrization of the solutions of the ARE in terms of the invariant subspaces of a certain matrix. The discrete-time counterpart of this parametrization was established in \cite{Wimmer-06-1}.
Reducing the problem to the analysis of a HARE is also one of the key ideas behind many theoretical results, see \cite{Clements-Wimmer,scherer,Wimmer-92,Wimmer-96} with  control applications including stochastic realization \cite{Picci-P-94,Ferrante-94-ieee,LPBook}, spectral facorization \cite{Ferrante-M-P-93,Ferrante-05-scl}, and smoothing \cite{Ferrante-P-98,Pavon-Wimmer}.
The interest for this topic continues in recent literature \cite{SANANDAMITADILIP20151,DILIP2017184}.
The advantage of considering homogeneous ARE is that it is possible to obtain a geometric picture describing a family of solutions. This family of solutions is parametrized in terms of invariant subspaces of a certain matrix.
Under specific assumptions this family is indeed the set of all the solutions of the HARE.
In particular, this is true when the reference solution of the ARE is stabilizing (or anti-stabilizing) so that the dynamics of the associate HARE is stable (or anti-stable). 
It is possible, however, to generalize this property to a reference solution that is {\em unmixing} i.e. such that the associated closed-loop matrix does not feature pairs of reciprocal eigenvalues.
In fact, this is the standing assumptions of most of the literature analysing the set of solutions of HARE's: the only exception to our knowledge is in \cite{Ferrante-Picci-TAC-17} where, however,  only the ARE associated to all-pass functions is considered.

In this paper we consider the following general HARE 
\beq\label{RicQ}
Q = A\tp Q A - A\tp Q B (R+B\tp QB)^{-1} B\tp QA
\eeq
where $A\in{\mathbb R}^{n\times n}$, $B\in {\mathbb R}^{n\times m}$ and $R=R\tp\in {\mathbb R}^{m\times m}$
and consider the general case where $A$ can be  {\em mixed} so that it can have pairs of eigenvalues 
$\lambda_1$ and $\lambda_2$ such that $\lambda_1 \lambda_2=1$.
Our contribution is threefold: first we show that in general
the HARE may have infinitely many families of solutions: each family is associated to one non-singular solution of \eqref{RicQ} and the solutions in any fixed family are parametrized in terms of the invariant subspaces of $A$. Also, we parametrize these families in terms of a linear equation. Second, we provide an explicit formula for the computation of the solutions of each family.
This formula is very simple and it proves to be useful even in the case when the ``unmixing'' assumption holds.
Third, we derive sufficient conditions under which the union of the families  covers the the whole set of solutions of \eqref{RicQ}.

\textbf{Notation.}  Given a matrix $M$, $M\tp$  denotes its transpose and $M\pinv$ its
Moore-Penrose pseudo-inverse.  The  kernel of  $M$ is denoted by  $\ker(M)$.

\section{Solutions of homogeneous Algebraic Riccati Equations}

Given the HARE \eqref{RicQ}, we only consider symmetric solutions so that when we say that $Q$ is a  solution of \eqref{RicQ} we mean that it is a symmetric matrix solving \eqref{RicQ}.
Let $\mathscr{I}_A$ be the set of $A$-invariant subspaces.
The following well-known (see \cite{Wimmer-06-1}) classical result parametrizes the set of solutions of 
\eqref{RicQ} in the case when $A$ is unmixed.
\begin{teor}\label{classic-t}
Let $(A,B)$ be a reachable pair and assume that $A$ is non-singular and that
$R=R\tp>0$. If $A$ is unmixed then
there is a bijective correspondence between the set of  solutions of
(\ref{RicQ}) and the set $\mathscr{I}_A$ of $A$-invariant subspaces. Such correspondence
is  defined by the map assigning to each solution $Q$ the $A$-invariant subspace $\ker(Q)$.
\end{teor}

In the following we relax the unmixing assumption.
As a first step we characterise the existence of invertible solutions of (\ref{RicQ}) and parametrize the set of such solutions in terms of a linear equation.
Consider the Stein (discrete-time Lyapunov) equation 
\beq\label{stein}
A P A\tp - P =BR^{-1}B\tp.
\eeq

\begin{lemm}\label{lemma-RL}
Assume that $A$ and $R$ are non-singular. There is a bijective correspondence between the set of non-singular  solutions of (\ref{RicQ}) and the set of non-singular solutions of \eqref{stein}. Such correspondence
is  defined by the map assigning to each non-singular solution $P$ of \eqref{stein}
the matrix $Q:=P^{-1}$ which is a non-singular  solution of (\ref{RicQ}).
\end{lemm}
\proof
Clearly, $Q_0$ is a non-singular solution of (\ref{RicQ}) if and only if
$$
Q_0^{-1} = A^{-1} [Q_0  -  Q_0 B (R+B\tp Q_0B)^{-1} B\tp Q_0]^{-1} A^{-\top}
$$
which, in view of the Sherman-Morrison-Woodbury matrix inversion formula, is equivalent to
$$
Q_0^{-1} = A^{-1} [Q_0^{-1}  + B R^{-1} B\tp ] A^{-\top}.
$$
The latter is clearly equivalent to $Q_0^{-1}$ being a non-singular solution of \eqref{stein}.
\qed

\begin{corr}
Let $(A,B)$ be a reachable pair. Then,
all the  solutions of equation \eqref{stein} are invertible so that given any solution
$P_0=P_0\tp$, the set $\mathscr{P}$ of all non-singular   solutions of \eqref{stein} is parametrized by 
$$
\mathscr{P}=\{P_\Delta:=P_0+\Delta:\ \Delta=\Delta\tp {\rm \ solves\ }A \Delta A\tp =\Delta\}.
$$
\end{corr}
\proof
Let $B_R:=BR^{-1/2}$. Since $(A,B)$ is reachable,  $(A,B_R)$ is reachable as well
so that in view of \cite[Lemma 3.1]{Ferrante-Ntog-Automatica-13}, all the  solutions of equation \eqref{stein} are invertible.
\qed

The following result shows that even when we drop all the assumptions, the kernel of any solution of \eqref{RicQ} is still an $A$-invariant subspace.

\begin{lemm}\label{lemma-ker-inv}
Let $Q$ be a   solution of (\ref{RicQ}).
Then $\ker(Q)\in \mathscr{I}_A$.
\end{lemm}
\proof 
Let $Q$ be fixed. We write (\ref{RicQ}) (which is now an identity) as
$Q= LQ(A-\varepsilon I)+\varepsilon LQ$ with $L:=A\tp  - A\tp Q B (R+B\tp QB)^{-1} B\tp$ 
and $\varepsilon$ being a constant such that both $A-\varepsilon I$ and
$I-\varepsilon L$ are non-singular.
This may be rewritten as $(I-\varepsilon L)Q(A-\varepsilon I)^{-1}=LQ$,
so that we immediately see that if $v\in\ker(Q)$ then $Q(A-\varepsilon I)^{-1}v=0$,
and hence $\ker(Q)$ is $(A-\varepsilon I)^{-1}$-invariant.
Thus $\ker(Q)$ is also $(A-\varepsilon I)$-invariant and, eventually,
$A$-invariant.
\qed

To any invertible solution $P_\Delta$ of \eqref{stein} we can associate a family of solutions of \eqref{RicQ}. This family is parametrized with respect to  the set $\mathscr{I}_A$  of $A$-invariant subspaces and in terms of the matrix $\Pi_\mathscr{S}$ that orthogonally projects into
$\mathscr{S}$. In formal terms, we have
 
\begin{teor}\label{famiglie}
Let $(A,B)$ be a reachable pair and assume that $A$ is non-singular and that
$R=R\tp>0$.
Let $P_\Delta$ be a solution of \eqref{stein}; let $\mathscr{S}\in\mathscr{I}_A$ and
$\Pi_\mathscr{S}$ be the orthogonal projector into $\mathscr{S}$.
Then
\beq\label{formulaesplperq}
Q=[(I-\Pi_\mathscr{S})P_\Delta(I-\Pi_\mathscr{S})]\pinv
\eeq
is a solution of  \eqref{RicQ}. Thus, for each given $P_\Delta\in\mathscr{P}$,
$$
\mathscr{Q}_\Delta:=\{Q=(I-\Pi_\mathscr{S})P_\Delta(I-\Pi_\mathscr{S}):\ \mathscr{S}\in\mathscr{I}_A \}
$$
defines a family of solutions of  \eqref{RicQ} parametrized in $\mathscr{I}_A$.
\end{teor}
\proof
Let $S$ be a matrix whose columns form an orthonormal basis for $\mathscr{S}$ and
$S_{\perp}$ be a matrix whose columns form an orthonormal basis for $\mathscr{S}^\perp$.
Consider a change of basis induced by the orthogonal matrix 
$T:=[S_{\perp}\mid S]$.
Clearly, $\Pi_\mathscr{S}=SS\tp$  and $\bar{A}:=T^{-1} AT$ has the form
$
\bar{A}=\bsmat A_1 & 0\\A_{21}&A_2\esmat.$
We partition $\bar{B}:=T^{-1} B$ conformably as $\bar{B}=\bsmat B_1 \\B_2\esmat$.
Finally we set
$\bar{\Pi}_\mathscr{S}:=T\tp \Pi_\mathscr{S} T=\bsmat 0 & 0\\0&I\esmat$ and
$\bar{P}_\Delta :=T\tp P_\Delta T$ which we partition conformably as
$
\bar{P}_\Delta=\bsmat P_1 & P_{12}\\P_{12}\tp &P_2\esmat.
$
Notice that since $(A,B)$ is reachable,  $(\bar{A},\bar{B})$ and hence
$(A_1,B_1)$ are also reachable.
Taking into account that $P_\Delta$ is by assumption a solution of \eqref{stein},
by a change of basis we immediately get 
\beq\label{steinbar}
\bar{A} \bar{P}_\Delta \bar{A}\tp - \bar{P}_\Delta =\bar{B}R^{-1}\bar{B}\tp.
\eeq
which, by employing the partitions just defined reads:
\beq
\label{??}
\bsmat A_1 & 0\\A_{21}&A_2\esmat \!\! \bsmat P_1 & P_{12}\\P_{12}\tp &P_2\esmat
 \!\! \bsmat A_1\tp & A_{21}\tp\\0&A_2\tp\esmat
-\bsmat P_1 & P_{12}\\P_{12}\tp &P_2\esmat  \!\! =  \!\!
\bsmat  {B_1}\\ {B_2}\esmat {\small R^{-1}}
\bsmat{B_1\tp} & {B_2\tp}\esmat.
\eeq
The upper-left block of this equation provides the following reduced-order Stein Equation
\beq\label{steinbar11}
A_1 P_1 A_1\tp - P_1=B_1 R^{-1}B_1\tp.
\eeq
Since $(A_1,B_1)$ is reachable, $P_1$ is non-singular \cite[Lemma 3.1]{Ferrante-Ntog-Automatica-13}
so that, in view of Lemma \ref{lemma-RL}, we have
\beq\label{RicQ-red}
P_1^{-1} = A_1\tp P_1^{-1} A_1 - A_1\tp P_1^{-1} B_1 (R+B_1\tp P_1^{-1} B_1)^{-1} B_1\tp P_1^{-1} A_1.
\eeq
By direct inspection, we can check that this implies that
$\bar{Q}:=\bsmat P_1^{-1} & 0\\0 &0\esmat$ is a solution of the ARE
\beq\label{RicQ-new-bas}
\bar{Q} = \bar{A}\tp \bar{Q} \bar{A} - \bar{A}\tp \bar{Q} \bar{B} (R+\bar{B}\tp \bar{Q} \bar{B})^{-1} \bar{B}\tp \bar{Q} \bar{A}.
\eeq
Therefore,
$$
Q:=T\bar{Q}T\tp$$
is a solution of \eqref{RicQ}.
It remains only to show that
$
Q=[(I-\Pi_\mathscr{S})P_\Delta (I-\Pi_\mathscr{S})]\pinv
$,
or equivalently that
$
Q\pinv=(I-\Pi_\mathscr{S})P_\Delta (I-\Pi_\mathscr{S}).
$
Since $T$ is orthogonal, we have 
$T\tp Q\pinv T=(T\tp Q T)\pinv$ so that it is sufficient to show that
\beq\label{eqfinalettqtpinv}
(T\tp Q T)\pinv=T\tp (I-\Pi_\mathscr{S})P_\Delta (I-\Pi_\mathscr{S})T
\eeq
The left-hand side of (\ref{eqfinalettqtpinv})
is
\beq\label{rhsef}
(T\tp Q T)\pinv=\bar{Q}\pinv=\bmat{cc}P_1^{-1} & 0\\0 &0\emat\pinv=
\bmat{cc}P_1 & 0\\0 &0\emat
\eeq
The right-hand side of (\ref{eqfinalettqtpinv})
is
\bea
\nn
&& T\tp (I-\Pi_\mathscr{S})P_\Delta (I-\Pi_\mathscr{S})T\\
\nn
&&= 
T\tp (I-\Pi_\mathscr{S})TT\tp P_\Delta TT\tp (I-\Pi_\mathscr{S})T\\
\nn
&&= \bmat{cc}I & 0\\0&0\emat\bmat{cc}P_1 & P_{12}\\P_{12}\tp &P_2\emat
\bmat{cc}I & 0\\0&0\emat=\bmat{cc}P_1 & 0\\0 &0\emat
\eea
which together with (\ref{rhsef}) proves (\ref{eqfinalettqtpinv}).
\qed

\begin{rem}
Notice that it is immediate to compute the kernel of $Q$ given by (\ref{formulaesplperq}):
$\ker(Q)=\ker([(I-\Pi_\mathscr{S})P_\Delta(I-\Pi_\mathscr{S})]\pinv)=\mathscr{S}$.
Therefore, when $A$ is unmixed and hence \eqref{stein} has exactly one solution $P_0$, in view of Theorem \ref{classic-t}, the only family 
\beq\label{par-figona}
\mathscr{Q}_0:=\{Q=(I-\Pi_\mathscr{S})P_0(I-\Pi_\mathscr{S}):\ \mathscr{S}\in\mathscr{I}_A \}
\eeq
provides the set of all the solutions of \eqref{RicQ}.
Even in this case, Theorem \ref{famiglie} is an important improvement with respect to Theorem \ref{classic-t} because of the explicit parametrization (\ref{par-figona})
that allows to compute
$
Q=[(I-\Pi_\mathscr{S})P_0(I-\Pi_\mathscr{S})]\pinv
$ 
as a function of the corresponding $A$-invariant subspace $\mathscr{S}$. 
\end{rem}

\begin{rem}
Assume that \eqref{stein} has non-singular solutions.
It is natural to ask whether or not the union of the solutions described in Theorem \ref{famiglie} covers the whole set of the solutions of \eqref{RicQ}.
Somehow surprisingly, the answer is negative as proven by the following counter-example.
Let $A:= 
\bsmat
 \frac{1}{2} & 0 & 0 & 0 \\
 0 & 2 & 0 & 0 \\
 1 & 0 & \frac{1}{2} & 0 \\
 0 & 1 & 0 & 2
\esmat$ and $B=I_4$.
By direct computation, we easily see that 
$Q_0:=\bsmat
 -\frac{87}{100} & 0 & \frac{9}{50} & 0 \\
 0 & \frac{21}{5} & 0 & \frac{9}{5} \\
 \frac{9}{50} & 0 & -\frac{27}{100} & 0 \\
 0 & \frac{9}{5} & 0 & \frac{27}{10}
\esmat$ is a non-singular solution of \eqref{RicQ}.
Therefore, our assumptions are satisfied.
On the other hand, we easily see that also
$Q_1:=
\bsmat
 -\frac{3}{13} & \frac{9}{13} & 0 & 0 \\
 \frac{9}{13} & \frac{12}{13} & 0 & 0 \\
 0 & 0 & 0 & 0 \\
 0 & 0 & 0 & 0
\esmat$ is  a solution of \eqref{RicQ}.
We now show, however, that such a solution is somehow {\em spurious} in the sense that it does not have the form in the right-hand side of \eqref{formulaesplperq}.
In fact, $Q_1\pinv=
\bsmat
 -\frac{4}{3} & 1 & 0 & 0 \\
 1 & \frac{1}{3} & 0 & 0 \\
 0 & 0 & 0 & 0 \\
 0 & 0 & 0 & 0
\esmat$ so that if $Q_1$ had the form 
of the right-hand side of \eqref{formulaesplperq}, then the corresponding matrix $P_\Delta$
should have the form 
$
P_\Delta=\bsmat
 -\frac{4}{3} & 1 & p_{13} & p_{14} \\
 1 & \frac{1}{3} & p_{23} & p_{24} \\
 p_{13} & p_{23} & p_{33} & p_{34}\\
 p_{14} & p_{24} & p_{34} & p_{44}
\esmat$ for suitable values of the entries $p_{ij}$.
It is, however, easy to check by direct inspection that with such a $P_\Delta$
the entry in position $(2,3)$ of the matrix $A P_\Delta A\tp - P_\Delta -BB\tp$ is equal to $2$ for every choice of 
the parameters $p_{ij}$. Therefore, $P_\Delta$ cannot be a solution of \eqref{stein}.
\end{rem}

The following result provides a sufficient condition ensuring that the  union of the solutions described in Theorem \ref{famiglie} covers the whole set of the solutions of \eqref{RicQ}.
\begin{teor}\label{unione-famiglie}
Let $(A,B)$ be a reachable pair and assume that $A$ is non-singular and that
$R=R\tp>0$. Moreover, assume that \eqref{stein} admits solutions. If $A$ has at most one pair of reciprocal eigenvalues and  these  eigenvalues (when present) have algebraic multiplicity equal to $1$ then each solution of \eqref{RicQ} is given by
\eqref{formulaesplperq} for a suitable solution $P_\Delta$    of \eqref{stein} and a suitable $\mathscr{S}\in\mathscr{I}_A$.
\end{teor}
\proof
Let $Q$ be a solution of \eqref{RicQ}.
In view of Lemma \ref{lemma-ker-inv}, we know that  $\ker(Q)$ is $A$-invariant. Then,
since we can perform the  change of basis described in the proof of Theorem \ref{famiglie},
we can assume, without loss of generality, that 
$
A=\bsmat A_1 & 0\\A_{21}&A_2\esmat$,
$B=\bsmat B_1 \\B_2\esmat$ and 
$Q=\bsmat Q_1 & 0\\0&0\esmat$ with $Q_1$ being non-singular.
By direct inspection, we can check that $Q_1$ is a non-singular solution
of the reduced-order ARE \eqref{RicQ-red} so that $P_1:=Q_1^{-1}$ is a solution of the 
reduced-order Stein equation \eqref{steinbar11}.
We only need to show that $P_1$ can be ``extended'' to a solution of \eqref{stein}, i.e. that  there exist matrices $X_{12}$ and $X_2=X_2\tp$ of suitable dimensions such that
$P:=\bsmat P_1 & X_{12}\\X_{12}\tp &X_2\esmat$ solves \eqref{stein}.
If $A_1$ is unmixed, equation \eqref{steinbar11} admits a unique solution; moreover,
 for any solution $P_0$ of  \eqref{stein} it is easy to check that the upper-left block of $P_0$ must satisfy \eqref{steinbar11} so that it is necessarily equal to $P_1$.
Consider now the case when $A_1$ is not unmixed. In this case, since $A$ has at most one pair of reciprocal eigenvalues which are both simple, not only is $A_2$ unmixed but we also have
\beq\label{crossunmixing}
\sigma(A_1)\cap  \sigma(A_2^{-1})=\emptyset.
\eeq
Now we consider \eqref{stein} with $P:=\bsmat P_1 & X_{12}\\X_{12}\tp &X_2\esmat$  as an equation in the unknowns $X_{12}$ and $X_2=X_2\tp$. If we write this equation block by block,  for the upper-left block we get equation \eqref{steinbar11} which is an identity,
for the upper-right block we get $A_1P_1A_{12}\tp+ A_1X_{12}A_{2}\tp-X_{12}=B_1R^{-1}B_2\tp$ which admits a solution $X_{12}$ because of \eqref{crossunmixing}.
Let $\bar{X}_{12}$ be such a solution and consider the upper-right block that now reads
$A_{12}P_1A_{12}\tp + A_2 \bar{X}_{12}\tp A_{12}\tp+A_{12}\bar{X}_{12}A_{2}\tp+
A_2X_2A_{2}\tp-X_2=B_2R^{-1}B_2\tp$ which admits a solution $X_{2}$ because
$A_2$ is unmixed.
\qed

\bibliographystyle{plain}

\end{document}